\documentclass[11pt]{article}

\usepackage{amsmath,amsfonts,amscd,amssymb,amsthm,enumitem,mathrsfs}

\usepackage[colorlinks=true,linkcolor=blue,citecolor=magenta]{hyperref}

\usepackage[justification=centering]{caption}

\newcounter{thmcount}[section]
\def\thethmcount{\thesection.\arabic{thmcount}}
\newtheorem{theorem}[thmcount]{Theorem}

\newtheorem{lemma}[thmcount]{Lemma}

\newtheorem{conjecture}[thmcount]{Conjecture}

\newtheorem*{conjecture*}{Conjecture}

\theoremstyle{definition}

\newtheorem{remark}[thmcount]{Remark}

\newtheorem{example}[thmcount]{Example}
\newtheorem*{remark*}{Remark}
\newtheorem*{problem*}{Problem}

\newcommand{\lra}{\longrightarrow}
\newcommand{\Zp}{{\mathbb{Z}_p}}
\newcommand{\Qp}{{\mathbb{Q}_p}}
\newcommand{\ZpG}{\mathbb{Z}_p[G]}

\newcommand{\CpG}{\mathbb{C}_p[G]}
\newcommand{\QpG}{\mathbb{Q}_p[G]}
\newcommand{\Cp}{{\mathbb{C}_p}}
\newcommand{\Ze}{{\mathbb{Z}}}
\newcommand{\Ce}{{\mathbb{C}}}
\renewcommand{\Re}{{\mathbb{R}}}

\newcommand{\PHj}{P_{(H,j)}}
\newcommand{\PtHj}{P^t_{(H,j)}}
\newcommand{\PJj}{P_{(J,j)}}
\newcommand{\PtJj}{P^t_{(J,j)}}
\newcommand{\Ptj}{P_{(t,j)}}
\newcommand{\Pttj}{P^t_{(t,j)}}
\newcommand{\Pntj}{P^t_{(n,j)}}
\newcommand{\Puk}{P_{(u,k)}}
\newcommand{\Pnk}{P_{(n,k)}}
\newcommand{\Qu}{\mathbb{Q}}
\newcommand{\QG}{\Qu[G]}
\newcommand{\calP}{\mathcal{P}}
\newcommand{\NT}{\mathrm{NT}}

\newcommand{\etncp}{{eTNC$_p$} }

\newcommand{\tensorZ}{\otimes_\Ze}
\newcommand{\tensorR}{\otimes_\Re}
\newcommand{\tensorZp}{\otimes_\Zp}
\newcommand{\tensorZl}{\otimes_\Zl}
\newcommand{\tensorQ}{\otimes_\Qu}
\renewcommand{\projlim}{{\lim\limits_{\longleftarrow}}}
\newcommand{\directlim}{{\lim\limits_{\longrightarrow}}}

\newcommand{\OkSp}{{\calO_{k, S}\left[\frac 1 p\right]}}
\newcommand{\TpFA}{{T_{p, F}(A)}}
\newcommand{\Cc}{{ C_{A,F}^{c, \bullet} }}
\newcommand{\Cf}{{ C_{A,F}^{f, \bullet} }}
\newcommand{\Cloc}{{ C_{A,F}^{loc, \bullet} }}
\newcommand{\lBK}{{ \lambda_{A,F}^{\mathrm{BK}} }}
\newcommand{\lNT}{{ \lambda_{A,F}^{\mathrm{NT}} }}
\newcommand{\lexp}{{ \lambda_{A,F}^{\mathrm{exp},j} }}
\newcommand{\sseq}{\subseteq}
\newcommand{\trivchar}{{{\bf 1}_G}}

\newcommand{\EG}{{\mathbb{E}[G]}}
\newcommand{\RG}{{\mathbb{R}[G]}}
\newcommand{\dett}{{\mathrm{det}}}

\renewcommand{\mod}{\mathrm{mod}\ }

\newcommand{\bd}{\mathbb D}
\newcommand{\bg}{\mathbb G}
\newcommand{\bh}{\mathbb H}
\newcommand{\bv}{\mathbb V}
\newcommand{\bzh}{\mathbb Z[\frac{1}{2}]}
\newcommand{\qbar}{\bar{\mathbb Q}}
\newcommand{\zbar}{\bar{\mathbb Z}}
\newcommand{\cf}{\mathcal F}
\newcommand{\ce}{\mathcal E}
\newcommand{\ca}{\mathcal A}

\newcommand{\A}{\mathfrak A}
\newcommand{\B}{\mathfrak B}

\newcommand{\D}{\mathrm D}
\newcommand{\cm}{\mathcal M}
\newcommand{\Up}{\mathcal O_{K,S_p}}

\newcommand{\calC}{\mathcal C}

\DeclareMathOperator{\diag}{diag}
\DeclareMathOperator{\Fit}{Fit}
\DeclareMathOperator{\Fr}{Fr}
\DeclareMathOperator{\Indshort}{I}
\DeclareMathOperator{\val}{val}
\DeclareMathOperator{\Nrd}{Nrd}
\DeclareMathOperator{\Ir}{Ir}

\newcommand{\calO}{\mathcal{O}}

\newcommand{\frakp}{\mathfrak{p}}

\newcommand{\calLstar}{\mathcal{L}^{\raisebox{1pt}{$\scriptstyle\star$}}}

\newcommand{\Lstar}{L^{\raisebox{2pt}{$\scriptstyle\star$}}} 

\newcommand{\vir}[1]{{\bf [} #1{\bf ]}}
\newcommand{\eins}{\boldsymbol{1}}
\DeclareMathOperator{\Det}{Det}
\DeclareMathOperator{\Ext}{Ext}
\DeclareMathOperator{\Ho}{H}
\DeclareMathOperator{\End}{End}
\DeclareMathOperator{\Pic}{Pic}
\DeclareMathOperator{\Spec}{Spec}
\DeclareMathOperator{\N}{N}
\newcommand{\fit}{\mathrm{Fit}_{\mathcal{O}}}
\DeclareMathOperator{\ind}{ind}
\DeclareMathOperator{\res}{res}
\DeclareMathOperator{\im}{im}

\newcommand{\CC}{\mathbb{C}}
\newcommand{\bc}{\mathbb{C}}
\newcommand{\FF}{\mathbb{F}}
\newcommand{\GG}{\mathbb{G}}
\newcommand{\NN}{\mathbb{N}}
\newcommand{\QQ}{\mathbb{Q}}
\newcommand{\Q}{\mathbb{Q}}
\newcommand{\RR}{\mathbb{R}}
\newcommand{\br}{\mathbb{R}}
\newcommand{\ZZ}{\mathbb{Z}}
\newcommand{\Z}{\mathbb{Z}}
\newcommand{\wrp}{\mathfrak{p}}
\newcommand{\ZG}{\mathbb{Z}[G]}

\newcommand{\calD}{\mathcal{D}}
\newcommand{\calE}{\mathcal{E}}
\newcommand{\calL}{\mathcal{L}}
\newcommand{\calM}{\mathcal{M}}
\newcommand{\calT}{\mathcal{T}}
\newcommand{\G}{\mathcal{G}}
\newcommand{\M}{\mathcal{M}}
\newcommand{\calN}{\mathcal{N}}
\newcommand{\co}{\mathcal{O}}
\newcommand{\calV}{\mathcal{V}}

\newcommand{\frakm}{\mathfrak{m}}
\newcommand{\frp}{\mathfrak{p}}
\newcommand{\frakA}{\mathfrak{A}}

\newcommand{\arch}{\mathrm{arch}}
\DeclareMathOperator{\aug}{aug}
\DeclareMathOperator{\Br}{Br}
\DeclareMathOperator{\centre}{Z}
\newcommand{\cok}{\text{cok}}
\newcommand{\cone}{\mathrm{cone}}
\newcommand{\et}{\mathrm{et}}
\newcommand{\even}{\mathrm{ev}}
\newcommand{\finite}{\mathrm{finite}}
\newcommand{\glob}{\mathrm{glob}}
\newcommand{\Imag}{\mathbf{I}}
\DeclareMathOperator{\inv}{inv}
\newcommand{\Irr}{\mathrm{Irr}}
\newcommand{\loc}{\mathrm{loc}}
\newcommand{\nr}{\mathrm{Nrd}}
\newcommand{\odd}{\mathrm{od}}
\newcommand{\old}{\mathrm{old}}
\newcommand{\perf}{\mathrm{p}}
\newcommand{\quot}{\mathrm{q}}
\newcommand{\real}{\mathrm{Re}}
\newcommand{\Real}{\mathbf{R}}
\newcommand{\Fitt}{\mathrm{Fitt}}
\newcommand{\pr}{\mathrm{pr}}
\newcommand{\Ann}{\mathrm{Ann}}
\newcommand{\bz}{\mathbb{Z}}
\newcommand{\La}{\Lambda}
\newcommand{\bq}{\mathbb{Q}}
\def\bigcapp{\raise1ex\hbox{\rotatebox{180}{$\biguplus$}}}
 \def\bigcappd{\raise1ex\hbox{\rotatebox{180}{$\displaystyle\biguplus$}}}

\def\Q{{\mathbb Q}}
\def\F{{\mathbb F}}
\def\Ql{{\mathbb Q}_l}
\def\Qp{{{\mathbb Q}_p}}
\def\Qb{\overline{\mathbb Q}}
\def\Qlb{{\overline{\mathbb Q}_l}}
\def\Qpb{{\overline{\mathbb Q}_p}}
\def\Z{{\mathbb Z}}
\def\Zp{{\mathbb Z}_p}
\def\Zl{{{\mathbb Z}_l}}
\def\R{{\mathbb R}}
\def\C{{\mathbb C}}
\def\N{{\mathbb N}}
\def\bL{{\mathbb L}}

\def\cO{{\mathcal O}}
\def\O{{\mathcal O}}
\def\cC{{\mathcal C}}
\def\cK{{\mathfrak K}}
\def\cF{{\mathfrak F}}

\def\Eta{{\mathrm H}}

\def\cL{{\mathcal L}}

\def\sA{{\mathscr A}}
\def\sB{{\mathscr B}}
\def\sL{{\mathscr L}}

\def\p{{\mathfrak p}}
\def\q{{\mathfrak q}}
\def\f{{\mathfrak f}}
\def\fF{{\mathfrak F}}
\def\fK{{\mathfrak K}}
\def\fG{{\mathfrak G}}
\def\fH{{\mathfrak H}}
\def\fL{{\mathfrak L}}
\def\fg{{\mathfrak g}}
\def\fh{{\mathfrak h}}

\newcommand{\hG}{\widehat{G}}

\newcommand{\Disc}{\operatorname{Disc}}
\newcommand{\Tr}{\operatorname{Tr}}
\newcommand{\Norm}{\operatorname{Norm}}
\newcommand{\ord}{\operatorname{ord}}
\newcommand{\GL}{\operatorname{GL}}
\newcommand{\SL}{\operatorname{SL}}
\newcommand{\PGL}{\operatorname{PGL}}
\newcommand{\Hom}{\operatorname{Hom}}
\newcommand{\Ind}{\operatorname{Ind}}
\newcommand{\Res}{\operatorname{Res}}
\newcommand{\sign}{\operatorname{sign}}
\newcommand{\rk}{\operatorname{rk}}
\newcommand{\corank}{\operatorname{corank}}
\newcommand{\coker}{\operatorname{coker}}
\newcommand{\codim}{\operatorname{codim}}
\newcommand{\cyc}{\operatorname{cyc}}
\newcommand{\Reg}{\operatorname{Reg}}
\newcommand{\Gal}{\operatorname{Gal}}
\newcommand{\Sel}{\operatorname{Sel}}
\newcommand{\Frob}{\operatorname{Frob}}
\newcommand{\BSD}{\operatorname{BSD}}
\newcommand{\tors}{\operatorname{tors}}
\newcommand{\LT}{\operatorname{LT}}
\newcommand{\id}{\mathrm{id}}
\newcommand{\Aut}{\operatorname{Aut}}
\newcommand{\lar}{\longrightarrow}
\newcommand{\too}{\longrightarrow}
\newcommand{\eps}{\varepsilon}
\newcommand{\isomto}{\overset{\sim\,}{\longrightarrow}}
\newcommand{\bn}[1]{\textcolor{blue}{#1}}
\newcommand{\rn}[1]{\textcolor{red}{#1}}
\newcommand{\gn}[1]{\textcolor{green}{#1}}
\newcommand{\iso}{\simeq}
\newcommand{\notdiv}{\nmid}
\newcommand{\triv}{\mathbf{1}}
\newcommand{\vsp}{\vspace{0.75em}}
\newcommand{\hol}{\mathrm{hol}}
\newcommand{\sing}{\mathrm{sing}}
\renewcommand{\subset}{\subseteq}
\renewcommand{\c}{\mathfrak{c}}
\renewcommand{\Re}{\operatorname{Re}}
\newcommand{\Art}{\operatorname{Art}}
\newcommand{\Infl}{\operatorname{Infl}}
\newcommand{\Cl}{\mathrm{Cl}}
\newcommand{\Mor}{\operatorname{Mor}}
\newcommand{\calK}{\mathcal{K}}
\newcommand{\calF}{\mathcal{F}}

\input cyracc.def       
\font\twlcyr=wncyr10 at 12pt
\def\sha{\text{\twlcyr\cyracc{Sh}}}

\begin{document}

\title{Numerical Evidence for a refinement of Deligne's Period Conjecture for Jacobians of Curves}
\author{Robert Evans, Daniel Macias Castillo and Hanneke Wiersema}



\newcommand{\Addresses}{{
  \bigskip
  \footnotesize

  R.~Evans, \textsc{University of Chichester, Bognor Regis Campus, Upper Bognor Road, Bognor Regis, West Sussex, PO21 1HR}\par\nopagebreak
  \textit{E-mail address}: \texttt{rob.evans@chi.ac.uk}

  \medskip

  D.~Macias Castillo, \textsc{Departamento de Matem\'aticas, 
Universidad Aut\'onoma de Madrid, 28049 Madrid (Spain);
and Instituto de Ciencias Matem\'aticas, 28049 Madrid (Spain).}\par\nopagebreak
  \textit{E-mail address}: \texttt{daniel.macias@uam.es}

  \medskip

  H.~Wiersema, \textsc{Department of Pure Mathematics and Mathematical Statistics, Centre for Mathematical Sciences, Wilberforce Road, Cambridge
CB3 0WB}\par\nopagebreak
  \textit{E-mail address}: \texttt{hw600@cam.ac.uk}

}}

\maketitle
%

\begin{abstract}

Let $A/\QQ$ be a Jacobian variety and let $F$ be a totally real, tamely ramified, abelian number field. Given a character $\psi$ of $F/\QQ$, Deligne's Period Conjecture asserts the algebraicity of the suitably normalised value $\calL(A,\psi,1)$ at $z=1$ of the Hasse-Weil-Artin $L$-function of the $\psi$-twist of $A$. We formulate a conjecture regarding the integrality properties of the family of normalised $L$-values $(\calL(A,\psi,1))_{\psi}$, and its relation to the Tate-Shafarevich group of $A$ over $F$. We numerically investigate our conjecture through $p$-adic congruence relations between these values. 

\end{abstract}
\let\thefootnote\relax\footnotetext{MSC: 11G10, 11G40 (primary), 11G30, 11G35 (secondary)}
\section{Introduction}
 
Deligne's Period Conjecture predicts that certain families of special $L$-values, after normalisation by appropriate periods, become algebraic numbers and moreover satisfy a natural Galois equivariance property. In many settings of interest, this conjecture constitutes a rationality prediction for normalised equivariant $L$-values.

In recent years, there has been much interest in the formulation and study of integral refinements of Deligne's Period Conjecture for the equivariant $L$-values that are associated to the base change of an abelian variety through a Galois extension of number fields. 
We refer the reader to \cite{Bley1,Bley2,Bley3,BleyMC1,BleyMC3,rbsd,BurnsMCWuthrich,DEW,WW}. However, as far as we are aware, any theoretical or numerical evidence obtained for such refinements has been restricted to the case of elliptic curves. 

The main purpose of this note is to formulate a conjectural integral refinement of Deligne's Period Conjecture in the case of Jacobian varieties of curves of higher genus, and to investigate it numerically through $p$-adic congruence relations between normalised $L$-values of their Dirichlet twists.  

In the case of elliptic curves defined over $\QQ$, Deligne's Period Conjecture itself is now known to be valid thanks to existing modularity results 
(cf. Remark \ref{modularity} below). However, the Jacobian varieties that we study do not, in general, correspond to a classical modular form.

To be a little more precise, we now let $A/\QQ$ be a Jacobian variety and $F$ be a totally real, tamely ramified, abelian number field. 
Fix also a rational prime $p$. We will then identify simple and explicit conditions on $A$, $F$ and $p$ which we conjecture are sufficient for a canonical normalised equivariant $L$-value of $A/F$, belonging a priori to the complex group algebra $\CC[\Gal(F/\QQ)]$, to in fact be $p$-adically integral. Moreover, we predict that its (Galois) action annihilates the $p$-primary Tate-Shafarevich group of $A$ over $F$.

We will study our explicit conditions through extended examples. 
We will also provide extensive supporting numerical evidence for our conjecture for Jacobians $A$ of curves of genus 2, base changed through number fields $F$ of degree $p$. In these situations, we will also make the additional partial consequences of our integrality prediction fully explicit.


Before stating our conjecture (as Conjecture \ref{Q1} below), we will define the relevant normalised $L$-values and briefly recall the precise statement of Deligne's Period Conjecture.

\subsection{Deligne's Period Conjecture}

In this section we recall the statement of Deligne's Period Conjecture.

We let $A/\QQ$ be an abelian variety of dimension $d$. Let $F$ be a totally real, tamely ramified, abelian number field, with Galois group $G:=\Gal(F/\QQ)$ and character group $\widehat{G}:=\Hom(G,\CC^\times)$. We denote the conductor of $F$ by $\mathfrak{f}$ and we write $C$ for the set of (rational) prime divisors of $\mathfrak{f}$.

Throughout this note, we will assume that the $C$-truncated Hasse-Weil-Artin $L$-series 
$$L_C(A,\psi,z)\,:=\,\prod_{\ell\,\nmid\, \mathfrak{f}}P_\ell(A,\psi,\ell^{-z})^{-1}$$ of $A$ and $\psi$,
has an analytic continuation to $z=1$. Here for each prime number $\ell$ and each character $\psi\in\widehat{G}$, 
we have written $P_\ell(A,\psi,t)$
for the Euler factor at $\ell$ for $h^1(A)\otimes[\psi]$, as defined in \cite[Not. 15]{DEW}.

We also associate the Gauss sum
\[ \tau^\ast(\psi) \,:=\, \sum_{a\in(\ZZ/\mathfrak{f}\ZZ)^\times}\psi(a)\cdot\zeta_\mathfrak{f}^{a}\] to each $\psi\in\widehat{G}$,
with $\zeta_\mathfrak{f}:={\rm exp}(2\pi i/\mathfrak{f})$. By abuse of notation, here we have written $\psi(a)$ for the image under $\psi$ of the restriction to $F$ of the automorphism of $\QQ(\zeta_\mathfrak{f})$ given by $\zeta_\mathfrak{f}\mapsto\zeta_\mathfrak{f}^{a}$. We then set
$$\calL_C(A,\psi)\,:=\,\frac{L_C(A,\psi,1)\cdot\tau^\ast(\check\psi)^d}{\Omega_A^+},$$ 
where $\check\psi$ is the contragredient character of $\psi$ and $\Omega_A^+$ is the real period of $A$, as defined in \cite[Section 2.4.2]{RE}.

Deligne's Period Conjecture, as formulated in \cite{Deligne1979}, predicts for each $\psi\in\widehat{G}$ a containment
\begin{equation}\label{1}\calL_C(A,\psi)\,\in\,\QQ(\psi)\end{equation} and, in addition, the Galois-equivariance property \begin{equation}\label{2}\calL_C(A,\fg\circ\psi)\;=\;\fg\big(\calL_C(A,\psi)\big)\end{equation} for all $\fg \in \Gal(\Q(\psi)/\Q)$. Here $\QQ(\psi)$ is the number field generated by the values of $\psi$. We refer the reader to the PhD thesis \cite{RE} of the first named author for details of how to relate these properties to Deligne's original formulation.


\subsection{An integral refinement of Deligne's Period Conjecture}

In this section we formulate the main conjecture of this article and we discuss some supporting evidence.

For each $\psi\in\widehat{G}$, consider now the idempotent $$e_\psi:=\frac{1}{|G|}\sum_{g\in G}\psi(g^{-1})g$$ of $\CC[G]$ at $\psi$. 
It is then a straightforward exercise to verify that properties (\ref{1}) and (\ref{2}) are valid for every $\psi$ in $\widehat{G}$, if and only if the element
\begin{equation}\label{3}\Theta_{C}(A_F)\,:=\,\sum_{\psi\in\widehat{G}}\calL_C(A,\psi)\cdot e_\psi\end{equation} of $\CC[G]$ belongs to $\QQ[G]$ (see, for instance, \cite[Lem. 2.8]{Bley1}).

Our numerical computations lead us to formulate the following conjecture, as a possible refinement of Deligne's Period Conjecture for Jacobian varieties. 

\begin{conjecture}\label{Q1} Assume that $A$ is the Jacobian variety of a smooth, projective curve over $\Q$, and has good reduction at all prime divisors of $\mathfrak{f}$. Fix an odd prime number $p$ that does not divide the product $\mathfrak{f}\cdot|A(F)_{\rm tor}| 
$.

Then the element $\Theta_{C}(A_F)$ defined in (\ref{3}) \begin{itemize}\item[(i)]
belongs to $\ZZ_{(p)}[G]$, and \item[(ii)] 
its action annihilates the $p$-primary Tate-Shafarevich group $\sha(A_F)[p^\infty]$ of $A$ over $F$.\end{itemize}\end{conjecture}

In \S \ref{primedegree} we will focus on the case where $F$ has degree $p$ and prime conductor. In this setting, we will make fully explicit some partial consequences of the integrality prediction in Conjecture \ref{Q1} (i). 

\begin{remark}\label{ConjRemark} For fixed $A$, $F$ and $p$, Conjecture \ref{Q1} simultaneously constitutes a refinement of the $a=0$ case of Prediction 8.1 in \cite{rbsd}, and an extension of the $a=0$ case of Prediction 8.4 in loc. cit., which was dependent upon additional hypotheses on reduction types and ramification. Let us note in passing that the general framework of loc. cit. may be used to relate Conjecture \ref{Q1} to the equivariant refinement of the Tamagawa number conjecture of Bloch and Kato \cite{BlochKato} that was formulated by Burns and Flach in \cite{BurnsFlach}. Conjecture \ref{Q1} is thus also closely linked to main conjectures in Iwasawa theory, 
although we will not explicitly discuss these connections here.\end{remark}

\begin{remark} In fact, it would be possible to use the approach of \cite{rbsd} to formulate an extension of Conjecture \ref{Q1} to more general number fields $F$. However, if for instance the archimedean place of $\QQ$ has non-trivial decomposition subgroup in $G$, then one would have to normalise each term $\calL_C(A,\psi)$ by a different period of $A$, depending on the parity of $\psi$. Also, the presence of wildly ramified primes in $F/\QQ$ would require replacing each Gauss sum $\tau^\ast(\psi)$ by a suitable `modified global Galois-Gauss sum' of $\psi$. For the computational purposes of this note, and the sake of simplicity, we have thus elected to restrict attention to totally real, tamely ramified, abelian number fields $F$.

We also refer the reader to Remark \ref{failureRk} below for a discussion of the weaker integrality properties for $\Theta_C(A_F)$ that we expect to hold for Jacobian varieties $A$ that have bad reduction at prime divisors of $\mathfrak{f}$. \end{remark}

Although the purpose of this note is to provide numerical evidence for Conjecture \ref{Q1} in higher dimension, one may provide the following theoretical evidence in support of Conjecture \ref{Q1} for elliptic curves $A$.

\begin{theorem}\label{Bley} Let $A/\QQ$ be an elliptic curve for which $L(A/\QQ,1)\neq 0$. Then there exist infinitely many primes $p$, and for each such $p$, there exist infinitely many fields $F$, satisfying the hypotheses of Conjecture \ref{Q1}, and such that Conjecture \ref{Q1} is valid for $A$, $F$ and $p$.\end{theorem}

Theorem \ref{Bley} follows upon combining the result \cite[Cor. 1.4]{Bley3} of Bley with the results \cite[Thm. 6.5, Prop. A.1]{rbsd} of Burns and the second named author and with the approach used to prove Thm. 8.6 in loc. cit. For brevity, we omit the details of the proof.

\begin{remark} Although the approach to proving Theorem \ref{Bley} outlined above only produces extensions $F/\QQ$ of $p$-power degree, 
it is possible to replace the use of Bley's result by the approach of \cite[\S 11.2]{rbsd} to prove a version of Theorem \ref{Bley} in which the degree of each extension $F/\QQ$ is divisible by arbitrarily many primes. 

However, the approaches of both \cite{Bley3} and of \cite[\S 11.2]{rbsd} rely crucially on the theory of modular symbols. Therefore, although one may be able to extend Theorem \ref{Bley} to higher-dimensional abelian varieties that are modular (associated to a newform of weight 2), these methods cannot be applied to non-modular Jacobian varieties.\end{remark}

Integral refinements of Deligne's Period Conjecture similar to that of Conjecture \ref{Q1}, both for values of the form (\ref{3}) and for analogous elements constructed from derivatives of Hasse-Weil-Artin $L$-series, have 
also been numerically investigated in the articles \cite{Bley1,Bley2,BleyMC1,BleyMC3,BurnsMCWuthrich}. However, as alluded to above, the investigations in these articles were undertaken exclusively for elliptic curves, satisfying moreover stricter hypotheses on reduction types than those that are in place in Conjecture \ref{Q1}.

In this note we provide extensive numerical evidence in support of claim (i) of Conjecture \ref{Q1} for Jacobian varieties of curves of genus 2, without any additional restrictions on reduction types beyond the assumed good reduction at all prime divisors of $\mathfrak{f}$.

We recall again that, through the general framework of \cite{rbsd}, one may show that our numerical evidence also supports aspects of the equivariant refinement of the Tamagawa number conjecture of Bloch and Kato \cite{BlochKato} that was formulated by Burns and Flach in \cite{BurnsFlach}.

Before stating our numerical evidence, we will briefly explain the kind of $p$-adic congruence relations that are encoded in Conjecture \ref{Q1}, and the consequences they can entail for the vanishing of $p$-primary Tate-Shafarevich groups. We hope that this explanation will serve as a first step, towards a higher-dimensional exploration of the kind of ideas that were applied in \cite{DEW} to study the arithmetic of elliptic curves. This will be the content of \S 2.

Finally, \S 3 we dedicate to the numerical investigation of Conjecture \ref{Q1}.

\subsubsection*{Acknowledgements}  The first author would like to thank Vladimir Dokchitser for the initial suggestion and valuable subsequent discussions on the problem of providing numerical evidence for the congruences that are subject of the present article, and for comments on an earlier version of the article.

The first and second author are also very grateful to David Burns for many helpful discussions.

The second author also wishes to thank Werner Bley for several pertinent observations on a preliminary version of the article, and Christian Wuthrich for his interest in this project.

The second author acknowledges support for this article as part of Grants CEX2019-000904-S and PID2019-108936GB-C21 funded by MCIN/AEI/ 10.13039/501100011033.

During the carrying out of this work the third author was supported by the Herchel Smith Postdoctoral Fellowship Fund, and the Engineering and Physical Sciences Research Council (EPSRC) grant EP/W001683/1.

\section{Congruence relations and Tate-Shafarevich groups}
\subsection{The general case}

The following result explains certain explicit consequences encoded within Conjecture \ref{Q1}.

For any vector $(x_\psi)_{\psi\in\widehat{G}}$ of complex numbers and any $g\in G$, we set
$$S_g\bigl((x_\psi)_{\psi\in\widehat{G}}\bigr)\,:=\,\sum_{\psi\in\widehat{G}}\psi(g)\cdot x_\psi\,\in\,\CC.$$

\begin{lemma}\label{CL}Fix $A$, $F$ and $p$ as in Conjecture \ref{Q1}. Then the following claims are valid.\begin{itemize}\item[(i)] Claim (i) of Conjecture \ref{Q1} holds if and only if for each $g\in G$, the sum
$S_g\bigl((\calL_C(A,\psi))_\psi\bigr)$ belongs to $|G|\cdot\ZZ_{(p)}$. \item[(ii)] Assume that Conjecture \ref{Q1} holds, that
$${L(A,\psi,1)}\neq 0$$ 
for every $\psi\in\widehat{G}$ and that for each $g\in G$, the sum
$S_g\bigl((\calL_C(A,\psi)^{-1})_\psi\bigr)$ belongs to $|G|\cdot\ZZ_{(p)}$. Then $\sha(A_F)[p^\infty]$ vanishes.
\end{itemize}\end{lemma}
\begin{proof} 
Claim (i) holds because for any vector $(x_\psi)_{\psi\in\widehat{G}}$ of complex numbers, the sum
\begin{equation}\label{gsumx}\sum_{\psi\in\widehat{G}}x_\psi\cdot e_\psi=\sum_{g\in G}\left(|G|^{-1}\sum_{\psi\in\widehat{G}}\psi(g)\cdot x_\psi\right)\cdot g^{-1}=\sum_{g\in G}\left(|G|^{-1}S_g\bigl((x_\psi)_{\psi}\bigr)\right)\cdot g^{-1}\end{equation} belongs to $\ZZ_{(p)}[G]$ if and only if, for each $g\in G$, the sum $S_g\bigl((x_\psi)_{\psi}\bigr)$ belongs to $|G|\cdot\ZZ_{(p)}$.

To prove claim (ii), we first observe that the stated non-vanishing hypothesis implies that 
$$L_C(A,\psi,1)\neq 0$$ 
for every $\psi\in\widehat{G}$ (note that $L_C(A,\psi,1)$ and $L(A,\psi,1)$ only differ by finitely many non-zero factors).

Given this fact, it is enough to show that the inverse
$$\Theta_C(A_F)^{-1}=\sum_{\psi\in\widehat{G}}\calL_C(A,\psi)^{-1}\cdot e_\psi$$ of $\Theta_C(A_F)$ in $\CC[G]$, also belongs to $\ZZ_{(p)}[G]$. By (\ref{gsumx}), the latter condition is valid if and only if, for each $g\in G$, the sum
$S_g\bigl((\calL_C(A,\psi)^{-1})_\psi\bigr)$ belongs to $|G|\cdot\ZZ_{(p)}$.
\end{proof}

\begin{remark} The condition that $$L(A,\psi,1)\neq 0$$ for every $\psi\in\widehat{G}$ is widely expected to hold, whenever the group $A(F)$ is finite. 
\end{remark}

\subsection{Extensions of prime degree}
\label{primedegree}

Throughout the rest of this article, we assume given odd prime numbers $p$ and $q$ such that $$q\equiv 1 \pmod{p}.$$
We then let $F$ be the (totally real, tamely ramified) unique subfield of $\Q(\zeta_q)$ that has degree $p$ over $\QQ$.

The following result makes some of the partial consequences of the integrality prediction in Conjecture \ref{Q1} fully explicit in this case.

\begin{lemma}\label{genus2conj}  Let $A$ be the Jacobian variety of a smooth, projective, rational curve. 
Assume that $A$ has good reduction at $q$, that $L(A/\Q,1)\neq 0$ and that $A(\QQ)$ contains no point of order $p$.

If claim (i) of Conjecture \ref{Q1} is valid for $A$, $p$ and $F$, then $L_{\{q\}}(A/\QQ,1)/\Omega_A^+$ belongs to $\ZZ_{(p)}$ and, for each non-trivial character $\psi$ of $G$, the following claims are valid:
\begin{enumerate}

\item[1.] (i) $\cL_q(A,\psi)$ belongs to $\Z_{(p)}[\zeta_p]$.

(ii) $\cL_q(A,\fg \circ \psi)\;=\;\fg\big(\cL_q(A,\psi)\big)\;$ for all $\fg \in \Gal(\Q(\zeta_p)/\Q).$

\medskip

\item[2.] $\cL_q(A,\psi)\;\equiv\; (-1)^d \cdot
\frac{L_{\{q\}}(A/\QQ,1)}{\Omega_A^+} \,\,\,\,\,\,\left(\mod\,\,\,(1-\zeta_p)\ZZ_{(p)}[\zeta_p]\right)$.

\end{enumerate}

\end{lemma}

\begin{proof} 

We extend each character $\psi$ of $G$ to a map $\psi:\CC[G]\to\CC$.

By \cite[Lem. 2.8]{Bley1}, the element $\Theta_q(A_F)$ belongs to $\QQ[G]$ if and only if for each $\psi \in \hG$, the element $\cL_q(A,\psi)$ belongs to $\QQ(\zeta_p)$ and satisfies claim 1. (ii).

If in addition $\Theta_q(A_F)$ belongs to $\ZZ_{(p)}[G]$, as predicted by claim (i) of Conjecture \ref{Q1}, then for any $\psi\neq\triv$ one has
$$\cL_q(A,\psi)=\psi\bigl(\Theta_q(A_F)\bigr)\,\in\,\psi\bigl(\ZZ_{(p)}[G]\bigr)=\Z_{(p)}[\zeta_p],$$ as required to prove claim 1. (i).
After observing that $\tau^\ast(\triv)=-1$, the same argument applied to the trivial character implies that $L_{\{q\}}(A/\QQ,1)/\Omega_A^+$ belongs to $\ZZ_{(p)}$, as was claimed. 

Now, for any $X$ in $\ZZ_{(p)}[G]$ and any $\psi\neq\triv$, one has
$$\psi(X)\,\equiv\,\triv(X)\,\,\,\,\,\,\left(\textrm{mod  }\psi(I_{p})\right),$$
where $I_p$ denotes the augmentation ideal in $\ZZ_{(p)}[G]$. Since $$\psi(I_p)=(1-\zeta_p)\ZZ_{(p)}[\zeta_p],$$
the congruence claimed in 2. clearly follows from this fact, with $X$ taken to be $\Theta_q(A_F)$, and after observing again that $\tau^\ast(\triv)=-1$.
\end{proof}

\begin{remark} The claim that $L_{\{q\}}(A/\QQ,1)/\Omega_A^+$ belongs to $\ZZ_{(p)}$ would be a consequence of the $p$-component of the Birch and Swinnerton-Dyer Conjecture for $A/\QQ$ {(after observing that $P_q(A,q^{-1})$ belongs to $\ZZ_{(p)}$)}.\end{remark}

\begin{remark} For each $\psi \in \hG$, we have 
\[\cL_q(A,\psi) \; 
\,= \frac{L_{\{q\}}(A,\psi,1)\cdot\overline{\tau^*(\psi)}^d}{\Omega_A^+}.\]
Note that if $\psi$ is non-trivial, then $q$ is totally ramified in the Artin field of $\psi$. Therefore, if $A$ has good reduction at $q$, using the Néron-Ogg-Shafarevich criterion, we find
\[
P_q(A,\psi,1/q)=1,
\]
and thus also
\begin{equation}\label{untruncated}\cL_q(A,\psi) \; 
\,= \frac{L(A,\psi,1)\cdot\overline{\tau^*(\psi)}^d}{\Omega_A^+}.\end{equation}
\end{remark}

\begin{remark}
\label{modularity}
If $A$ is an elliptic curve, it  follows from results of Shimura \cite{shimura}, in combination with the modularity of $A$ \cite{wiles,tw,BCDT}, that Deligne's period conjecture is valid for each twist of $A$ by a Dirichlet character $\psi$. For an explicit statement of this claim, and a proof for some additional cases of Artin representations, see the article of Bouganis and Dokchitser \cite{bd}.

Moreover, in this case, one further knows that the element $\cL_q(A,\psi)$ belongs to $\Z[\zeta_p]$ for any non-trivial character $\psi$ of $G$ which validates Stevens's Manin constant conjecture. For details of this assertion, see \cite[Thm. 2 a)]{WW}.

We note that such a containment is slightly stronger than the containment of $\cL_q(A,\psi)$ in $\Z_{(p)}[\zeta_p]$ that is encoded in Conjecture \ref{Q1}, and that our numerical computations in \S \ref{numerical} below have confirmed this stronger integrality property in all appropriate examples of Jacobian varieties.

\end{remark}

\section{Numerical Evidence}
\label{numerical}

In this section we will provide both extended examples and data in support of the explicit predictions discussed in previous sections.

Before we provide the examples, let us note that all calculations were carried out in MAGMA \cite{magma} with the precision set to at least $10$, and that all (untruncated) $L$-values were computed using the algorithm of T. Dokchitser \cite{Lvaluecomp}. To determine the $L$-value, note that $\Z[\zeta_p]$ can be identified with a \textit{discrete} subgroup of $\C^{p-1}$ and in this way we can test whether the modified $L$-value is close to a point in the image of $\Z[\zeta_p]$.

We also stress that many of these computations are subject to our assumption that the corresponding $L$-series admit an analytic continuation to $z=1$. 

To numerically verify the full extent of Conjecture \ref{Q1} for examples of triples $(A,F,p)$ for which the $p$-primary Tate-Shafarevich group of $A$ over $F$ does not vanish, such as those discussed in \S \ref{nonvanishing} below, one would have to explicitly determine the Galois structure of this group. This seems, to us, to be a very delicate problem in any example. Nevertheless we will attempt to return to this problem in future work, and we expect that the data provided in the tables below will be helpful to carry out such full numerical verifications.

\subsection{Explicit examples}

We first provide an explicit illustration of claim (i) of Conjecture \ref{Q1}. 

\begin{example}\label{ex1}($p=5,\,q=11$)
Let $A$ be the Jacobian variety of the genus 2 curve given by the equation
\[y^2 + (x^3 + 1)y \;=\; x^5 - x^4 - 5x^3 + 4x^2 + 4x - 4,\]
with LMFDB label 427.a.2989.1. We find good reduction at $11$, and no non-trivial $5$-torsion points in $A(\QQ)$.

We let $F$ be the degree 5 subfield of $\Q(\zeta_{11})$ and let $\sigma$ denote the automorphism of $\Q(\zeta_{11})$ defined by $\sigma(\zeta_{11})=\zeta_{11}^2$. We abuse notation and denote by $\sigma$ its restriction to $F$. We then let $\psi_j$ be the linear character on $\Gal(F/\Q)$ mapping $\sigma$ to $\zeta_5^j.$ The Dirichlet character corresponding to $\psi_j$ is the unique group homomorphism $(\Z/11\Z)^\times \to \C^\times$ that maps $2$ to $\zeta_5^j$.

The modified, normalised $L$-values of $A$ are
\begin{align*}
\cL_{11}(A,\triv)&\;=\;P_{11}\Big(A,\frac{1}{11}\Big)\cdot(-1)^2\frac{L(A/\QQ,1)}{\Omega_A^+}\\[0.5em]
&\;=\;\left(1-\frac{1}{11}+\frac{2}{11^2}-\frac{11}{11^3}+\frac{11^2}{11^4}\right)\cdot 0.01020408163
\;\approx\;\frac{2^3}{7\cdot 11^2}=:\alpha_0,\\[0.5em]
\cL_{11}(A,\psi_1)& \;\approx\; -2(1+\zeta_5+\zeta_5^3)=:\alpha_1,\,\,\,\,\,\,\,\,\,
\cL_{11}(A,\psi_2) \;\approx\;  -2(1+\zeta_5+\zeta_5^2)=:\alpha_2, \\
\cL_{11}(A,\psi_3)& \;\approx\;  2(\zeta_5+\zeta_5^2)=:\alpha_3, \,\,\,\,\,\,\,\,\,\,\,\,\,\,\,\,\,\,\,\,\,\,\,\,
\cL_{11}(A,\psi_4) \;\approx\;  2(\zeta_5+\zeta_5^3)=:\alpha_4.
\end{align*}
where `$\approx$' means `equal to $10$ significant figures'. 

We find that $\Theta_{11}(A_F)$ belongs to $\Z_{(5)}[G]$, and thus that claim (i) of Conjecture \ref{Q1} is valid. Explicitly, 
\[
\Theta_{11}(A_F)= \frac{2}{7 \cdot 11^2}\left(-(2  \cdot 13^2) (1+\sigma+\sigma^4)+509 (\sigma^2+\sigma^3)\right).\
\]

Let us finally make explicit the congruences given in claim 2. of Lemma \ref{genus2conj}.
On the one hand we have
\begin{align*}
\alpha_1\;=\;-2(1+\zeta_5+\zeta_5^3)
&\;=\; 2(1-\zeta_5)(2+\zeta_5+\zeta_5^2)-6\\
&\;\equiv\; -1  \,\,\,\,\,\,\left(\mod\,\,\,(1-\zeta_5)\ZZ_{(5)}[\zeta_5]\right),
\end{align*}
whilst on the other hand, we have
\begin{align*}
\alpha_0 &\;\equiv\; \frac{2^3}{2\cdot 1^2} \,\,\,\,\,\,\left(\mod\,\,\,(1-\zeta_5)\ZZ_{(5)}[\zeta_5]\right) 
&\;\equiv\; -1 \,\,\,\,\,\,\left(\mod\,\,\,(1-\zeta_5)\ZZ_{(5)}[\zeta_5]\right).
\end{align*}
One readily computes that $\alpha_2,\alpha_3$ and $\alpha_4$ are also congruent to $-1$. 
\end{example}

\subsubsection{Failure of the integrality claim}
Let us now give some examples where the hypotheses of Conjecture \ref{Q1} fail to hold, and so does its integrality claim.

\begin{example}\label{ex2}($p=3,\,q=7$)
Let $A$ be the Jacobian variety of the {genus 2} curve given by the equation
\[
y^2 + (x^3 + 1)y = -2x^4 + 4x^2 - 9x - 14
\]
with LMFDB label 294.a.8232.1. Then $A$ has both bad reduction at $7$ and points of order 3 over $\QQ$.

Let $F$ be the degree 3 subfield of $\Q(\zeta_7)$ and let $\sigma$ be the automorphism of $\Q(\zeta_7)$ defined by $\sigma(\zeta_7)=\zeta_7^2$. We abuse notation and denote by $\sigma$ its restriction to $F$. Let $\psi$ to be the linear character of $F$ mapping $\sigma$ to $\zeta_3$. 

We compute
\[\mathcal{L}_7(A,\triv) \; = \; \frac{1}{7^2} \qquad \text{and} \qquad \mathcal{L}_7(A,\psi) \; = \; - (1+2\zeta_3).\]
We thus find 
\[
\Theta_7(A_F)=\frac{1}{3\cdot 7^2}  \left(1 - (2 \cdot 73) \cdot\sigma + (2^2 \cdot 37) \cdot\sigma^2\right) \not \in \Z_{(3)}[G].
\]

Let us also note that $\mathcal{L}_7(A,\psi)=-(1+2 \zeta_3)=-\zeta_3 (1-\zeta_3)$ but that $\calL_7(A,\triv)$ is not congruent to $0$ modulo $(1-\zeta_3)$, so the congruence in claim 2. of Lemma \ref{genus2conj} also fails to hold.
\end{example}

\begin{example}[$p=3,\,q=7$]
 \label{27a3}
We use the notation of Example \ref{ex2}. The elliptic curve $E/\Q$ with Cremona label 27a3 has good reduction at 7, but does have points of order $3$ over $\QQ$.

We compute
\[\mathcal{L}_7(E,\triv) \; = \; -\frac{1}{7} \qquad \text{and} \qquad \mathcal{L}_7(E,\psi) \; = \; 1.\]
We thus find 
\[
\Theta_7(E_F)=\frac{1}{3\cdot 7} \left(13 - 2^3\cdot \sigma - 2^3 \cdot \sigma^2\right) \not \in \Z_{(3)}[G].
\]
\end{example}

\begin{example}[$p=3,\,q=7$]
 \label{147b1}
We use the notation of Example \ref{ex2}. The elliptic curve $E/\Q$ with Cremona label 147b1 has no non-trivial $3$-torsion over $\QQ$, but does have bad reduction at $7$.

We compute
\[\mathcal{L}_7(E,\triv) \; = \; -1 \qquad \text{and} \qquad \mathcal{L}_7(E,\psi) \; = \; \frac{7}{13}.\]
We thus find
\[
\Theta_7(E_F)=\frac{1}{3\cdot 13}  \left(1  -20 \cdot \sigma -20 \cdot \sigma^2\right) \not \in \Z_{(3)}[G].
\]
\end{example}

\begin{remark}\label{failureRk} Even in the (non-equivariant) setting of the {Birch and Swinnerton-Dyer Conjecture}, the presence of $p$-torsion rational points leads to the failure of integrality at $p$ for normalised $L$-values.

Let $A/\QQ$ be a Jacobian variety such that $A(\QQ)[p]=0$ but which has bad reduction at the ramified place $q$. Then the classical periods and Gauss sums occurring in the definition of $\Theta_{q}(A_F)$ can differ from the determinants of canonical `$K$-theoretical periods' that occur naturally in the formulation of refined conjectures of Birch and Swinnerton-Dyer type (see \cite{rbsd}). In our specific setting, such a difference is bounded by a factor of $p^d$ (times a rational $p$-adic unit).

In this case, it would therefore be natural to expect that, after 
multiplication by the additional factor $p^d$, the element $\Theta_{q}(A_F)$ would belong to $\ZZ_{(p)}[G]$.

Such an expectation is consistent with our numerical computations but, since it does not necessarily pin down the sharpest possible bound on the denominators at $p$ that can occur in $\Theta_{q}(A_F)$, we leave its precise formulation and its thorough numerical investigation to future work.

\end{remark}

\subsubsection{Non-vanishing $p$-primary Tate-Shafarevich groups}\label{nonvanishing} 
We finally provide some examples for which claim (i) of Conjecture does hold, but the $p$-primary Tate-Shafarevich group of $A$ over $F$ does not vanish. In all such instances, we have found that the condition on the sums $S_g\bigl((\calL_C(A,\psi)^{-1})_\psi\bigr)$ that is stated in Lemma \ref{CL} (ii) fails to hold. Therefore, the non-vanishing of $\sha(A_F)[p^\infty]$ does not provide a counterexample to claim (ii) in Conjecture \ref{Q1}.


\begin{example}($p=5,\,q=11$)
\label{11a2}
We use the notation of Example \ref{ex1}.
The elliptic curve $E/\QQ$ with Cremona label 1246b2 has good reduction at $11$ and no $5$-torsion over $\QQ$. The analytic order of $\sha(E_F)$ is $625=5^4$.

We compute
\begin{multline*}
\mathcal{L}_{11}(E,\triv)  \; = \; \frac{-300}{11},\\
\mathcal{L}_{11}(E,\psi_1)  \; = \;   -8 \zeta_5^3 - 22 \zeta_5^2 - 22 \zeta_5 - 8 ,\,\,\,\,\,\,\,\,\,\,\,\,
\mathcal{L}_{11}(E,\psi_2)  \; = \;  22\zeta_5^3 + 14\zeta_5 + 14  \\
\mathcal{L}_{11}(E,\psi_3)  \; = \; -14\zeta_5^3 + 8\zeta_5^2 - 14\zeta_5,\,\,\,\,\,\,\,\,\,\,\,\,\,\,\,\,\,\,\,\,\,\,\,\,
\mathcal{L}_{11}(E,\psi_4)  \; = \; 14\zeta_5^2 + 22\zeta_5 + 14   
\end{multline*}

We find that $\Theta_{11}(E_F)$ belongs to $\Z_{(5)}[G]$. Explicitly,
\[
\Theta_{11}(E_F)= \frac{2}{11}\left(-2^3\cdot (1+\sigma^3)+(5\cdot 17) \cdot (\sigma^1+\sigma^2)+ (2^2 \cdot 3^2)\cdot \sigma^4\right).
\]

However, the sum $S_g\bigl((\calL_{11}(E,\psi)^{-1})_\psi\bigr)$ does not belong to $5\cdot\ZZ_{(5)}$ for all $g \in G$. For example, taking $g=1$, we find
\[
S_1\bigl((\calL_{11}(E,\psi)^{-1})_\psi\bigr)=\frac{739}{300}.
\]
\end{example}

\begin{example}($p=3,\,q=13$)
The elliptic curve $E/\Q$ with Cremona label $448c5$ has good reduction at 13 and no non-trivial $3$-torsion over $\QQ$. The analytic order of $\sha(E_F)$ is $729=3^6$.

Let $F$ be the degree 3 subfield of $\Q(\zeta_{13})$, and let $\sigma$ be the automorphism of $\Q(\zeta_{13})$ defined by $\sigma(\zeta_{13})=\zeta_{13}^2$. Let $\psi$ be the linear character of $F$ mapping $\sigma$ to $\zeta_3$.

We compute
\[\mathcal{L}_{13}(E,\triv) \; = \; \frac{-90}{13} \qquad \text{and} \qquad \mathcal{L}_{13}(E,\psi) \; = \;- 9 \cdot \zeta_3.\]
We find that $\Theta_{13}(E_F)$ belongs to $\ZZ_{(3)}[G]$. Explicitly,
\[
\Theta_{13}(E_F)=\frac{1}{13}\left(9 \cdot 1 -108\cdot \sigma +9  \cdot \sigma^2\right).
\]

However, the sum $S_g\bigl((\calL_{13}(E,\psi)^{-1})_\psi\bigr)$ does not belong to $3\cdot\ZZ_{(3)}$ for all $g \in G$. For example, taking $g=1$, we find
\[
S_1\bigl((\calL_{13}(E,\psi)^{-1})_\psi\bigr)=\frac{-1}{30}.
\]
\end{example}

\subsection{Further numerical evidence}

In this section we will provide a guide to the tables in the following section.

Let $S$ be the set of pairs $(p,q)$ of odd primes
\[\{(3,7),(3,13),(3,19),(3,31),(5,11),(5,31),(7,29)\},\]
which all satisfy $q\equiv 1 \bmod{p},$ let $F_{p,q}$ be the degree $p$ subfield of $\Q(\zeta_q)$ and write $G=\Gal(F_{p,q}/\Q).$ For each of the 38 abelian varieties $A$ of conductor at most $500$ that arise as Jacobians of genus 2 curves over $\Q$ that are listed in the LMFDB \cite{LMFDB} and each pair $(p,q) \in S,$ we calculated the $p$-tuple of modified $L$-values 
\[\big(\cL_q(A,\psi) \; : \; \psi\in \hG\,\big).\]

For presentational reasons, we only list one member of the above tuple in the tables below; however, all $L$-values were computed independently of one another.  Moreover, by Lemma \ref{genus2conj} 1. (ii), if the congruence relation in claim 2. of Lemma \ref{genus2conj} holds for one non-trivial character $\psi \in \hG$, then it follows for \textit{all} non-trivial linear characters.

We remark that some of the curves listed below give rise to isogenous Jacobian varieties and so their (unmodified) $L$-values are equal. However, the hypothesis on the rational $p$-torsion subgroups in Lemma \ref{genus2conj} is \textit{not} isogeny-invariant. Indeed, we see from Table \ref{gen2tab2} below that, if $A$ (resp. $A'$) is the Jacobian variety of curve 277.a.277.1 (resp. 277.a.277.2), then $A$ and $A'$ are isogenous but $A(\Q)[3]$ has order 3 whereas $A'(\Q)[3]$ is trivial and, as the table shows, the congruence relation in claim 2. of Lemma \ref{genus2conj} is (numerically) false for $A$ and (numerically) true for $A'.$ 

\subsubsection{How to read the tables}
Note that we are using the same curve-labels as those used in the LMFDB \cite{LMFDB}. We write $N_A$ for the conductor of $A$ in each case.

The tables below are colour-coded so that Lemma \ref{genus2conj} can be verified `at a glance'. Columns 2, 3 and 7 represent the hypotheses $q\nmid N_A$, $p\nmid |A(\Q)_\mathrm{tors}|,$ and integrality of $\Theta_C$ respectively, and are coloured \textcolor{blue}{blue} if they are \textit{not} satisfied (with `N' for no in column 7). Columns 4, 5 and 6 represent the conclusions 
and are coloured \textcolor{red}{red} if they do \textit{not} hold. More precisely, we write $\p$ for the prime ideal $(1-\zeta_p)\ZZ[\zeta_p]$ of $\ZZ[\zeta_p]$ and $\ord_\p$ for the normalised $\p$-adic valuation on $\QQ(\zeta_p)$. Then
column 6 (headed `$\ord_\p$') denotes 
\[\ord_\p\left(\cL_q(A,\psi)-\frac{L_{\{q\}}(A/\QQ,1)}{\Omega_A^+}\right).\] 

In terms of this colour-coding, Lemma \ref{genus2conj} can be rephrased as follows:

\begin{enumerate}[leftmargin=3em]

\item[1)] \fbox{columns 2, 3 and 7 are black} $\;\Rightarrow\;$ \fbox{column 4 is black}\,,

\item[2.i)] \fbox{column 2, 3 and 7 are black} $\;\Rightarrow\;$ \fbox{column 5 is black}\,,

\item[3)] \fbox{columns 2, 3 and 7 are black} $\;\Rightarrow\;$ \fbox{column 6 is black}\,.
\end{enumerate}

Conversely, we can rephrase Conjecture \ref{Q1} i) as follows
\begin{itemize}

\item \fbox{columns 2 and 3 are black} $\;\Rightarrow\;$ \fbox{column 7 is black}\,.

\end{itemize}

\pagebreak

\subsubsection{Tables}\label{tables}
\begin{table}[h!]
\centering
\caption{$p=3,\,q=7,\,\psi(3)=\zeta_3,\,\p=(1-\zeta_3).$}
\scriptsize
\begin{tabular}{|c|c|c|c|c|c|c|} \hline

Curve Label      &              $N_A$       &  $|A(\Q)_{\mathrm{tors}}|$  &       $\cL(A,\triv)$      &              $\cL(A,\psi)$      &                                              $\ord_\p$      &              $\Theta_C$ \\ \hline \hline                                                                                                                           
169.a.169.1      &              $13^{2}$                  &   $19^{1}$                       			&                             $3^{1}7^{-2}19^{-1}$&              0                      				 &                                                                     2      &              \\ \hline
196.a.21952.1      &              $2^{2}{\color{blue}7^{2}}$  &   $2^{2}{\color{blue}3^{2}}$ &        ${\color{red}3^{-1}}7^{-2}$  &                           $1+\zeta_{3}$      &                                              \textcolor{red}{-2}      &              \textcolor{blue}{N} \\ \hline
249.a.249.1      &              $3^{1}83^{1}$             & $2^{1}7^{1}$                       		&                             $2^{1}7^{-3}$      &              $1+\zeta_{3}$      &                                               1      &              \\ \hline
249.a.6723.1      &              $3^{1}83^{1}$             & $2^{2}7^{1}$                       		&                              $2^{1}7^{-3}$      &              $1+\zeta_{3}$      &                                              1      &              \\ \hline
256.a.512.1      &              $2^{8}$                 &  $2^{2}5^{1}$                       		&                           $5^{-1}7^{-2}$      &              $1+\zeta_{3}$      &                                                1      &              \\ \hline
277.a.277.1      &              $277^{1}$               &  ${\color{blue}3^{1}}5^{1}$&   $5^{-1}7^{-2}$      &              $1+\zeta_{3}$      &                                                1      &              \\
\hline

\end{tabular}
\label{gen2tab1}
\end{table}

\pagebreak

\begin{table}[h!]
\centering
\caption{$p=3,\,q=13,\,\psi(2)=\zeta_3,\,\p=(1-\zeta_3).$}
\scriptsize
\begin{tabular}{|c|c|c|c|c|c|c|} \hline

Curve Label        &     $N_A$        &           $|A(\Q)_{\mathrm{tors}}|$  &    $\cL(A,\triv)$        &     $\cL(A,\psi)$        &          $\ord_\p$        &     $\Theta_C$ \\ \hline \hline
                                                                                                                         
169.a.169.1        &     ${\color{blue}13^{2}}$      &     $19^{1}$                       				  &                                $13^{-1}19^{-1}$        &     $\zeta_{3}$  &                            1        &     \\\hline
196.a.21952.1        &     $2^{2}7^{2}$               &    $2^{2}{\color{blue}3^{2}}$        &  $3^{1}13^{-2}$        &     $3$        &                    4        &     \\\hline
249.a.249.1        &     $3^{1}83^{1}$        &  $2^{1}7^{1}$        &       $2^{1}3^{1}7^{-1}13^{-2}$        &     0        &      2         &     \\\hline
249.a.6723.1        &     $3^{1}83^{1}$        &  $2^{2}7^{1}$        &       $2^{1}3^{1}7^{-1}13^{-2}$        &     0        &      2        &     \\\hline
256.a.512.1        &     $2^{8}$        &        $2^{2}5^{1}$        & $13^{-2}$        &     $\zeta_{3}$               &      			           1        &     \\\hline 
277.a.277.1        &     $277^{1}$        &       ${\color{blue}3^{1}}5^{1}$       &       $3^{1}5^{-1}13^{-2}$        &     $\zeta_{3}$        &            \textcolor{red}{0}        &     \textcolor{blue}{N} \\\hline
277.a.277.2        &     $277^{1}$        &        $5^{1}$                       			  &                    $3^{3}5^{-1}13^{-2}$        &     $9\zeta_{3}$        &           4        &     \\
\hline

\end{tabular}
\label{gen2tab2}
\end{table}

\pagebreak

\begin{table}[h!]
\centering
\caption{$p=3,\,q=19,\,\psi(2)=\zeta_3,\,\p=(1-\zeta_3).$}
\scriptsize
\begin{tabular}{|c|c|c|c|c|c|c|} \hline

Curve Label        &     $N_A$        &       $|A(\Q)_{\mathrm{tors}}|$ &       $\cL(A,\triv)$        &     $\cL(A,\psi)$        &            $\ord_\p$        &    $\Theta_C$ \\ \hline \hline
                                                                                                                            
169.a.169.1        &     $13^{2}$        &       $19^{1}$           & $3^{3}19^{-3}$        &     $3(1+\zeta_{3})$     &         2        &     \\\hline
196.a.21952.1        &     $2^{2}7^{2}$        & $2^{2}{\color{blue}3^{2}}$   &          $3^{1}19^{-2}$        &     $12\zeta_{3}$  &               3        &     \\\hline
249.a.249.1        &     $3^{1}83^{1}$        &    $2^{1}7^{1}$           &$2^{1}3^{2}7^{-1}19^{-2}$        &     $3(1+\zeta_{3})$ &                          2        &     \\\hline
249.a.6723.1        &     $3^{1}83^{1}$        &    $2^{2}7^{1}$           &$2^{1}3^{2}7^{-1}19^{-2}$        &     $3(1+\zeta_{3})$&                           2        &     \\\hline
256.a.512.1        &     $2^{8}$        &          $2^{2}5^{1}$               &   $2^{-1}5^{-1}13^{1}19^{-2}$        &     $\zeta_{3}$&             1        &     \\\hline
277.a.277.1        &     $277^{1}$        &    ${\color{blue}3^{1}}5^{1}$     &     $2^{3}5^{-1}19^{-2}$        &     $1$        &                  2        &     \\\hline
277.a.277.2        &     $277^{1}$        &      $5^{1}$                      &		                $2^{3}3^{2}5^{-1}19^{-2}$        &     $9$      &               6        &     \\
\hline

\end{tabular}
\label{gen2tab3}
\end{table}

\pagebreak

\begin{table}[h!]
\centering
\caption{$p=3,\,q=31,\,\psi(3)=\zeta_3,\,\p=(1-\zeta_3).$}
\scriptsize
\begin{tabular}{|c|c|c|c|c|c|c|} \hline

Curve Label        &     $N_A$        &                 $|A(\Q)_{\mathrm{tors}}|$  &      $\cL(A,\triv)$        &     $\cL(A,\psi)$        &    $\ord_\p$        &     $\Theta_C$ \\ \hline \hline
                                                                                                                         
169.a.169.1        &     $13^{2}$        &       $19^{1}$                         				&                       $2^{4}3^{1}19^{-1}31^{-2}$        &     0        &         2        &     \\\hline
196.a.21952.1        &     $2^{2}7^{2}$        &$2^{2}{\color{blue}3^{2}}$        & $2^{2}3^{1}31^{-2}$        &     0        &              2        &     \\\hline
249.a.249.1        &     $3^{1}83^{1}$        &    $2^{1}7^{1}$              &     $2^{5}7^{-1}31^{-2}$        &     $-1$        &              4        &     \\\hline
249.a.6723.1        &     $3^{1}83^{1}$        &    $2^{2}7^{1}$              &     $2^{5}7^{-1}31^{-2}$        &     $-1$        &             4        &     \\\hline
256.a.512.1        &     $2^{8}$           &			$2^{2}5^{1}$                  					&                  $2^{3}31^{-2}$        &     $2$        &     	                  2        &     \\\hline
277.a.277.1        &     $277^{1}$        &${\color{blue}3^{1}}5^{1}$                         	&                        $2^{1}3^{1}31^{-2}$        &     0        &               2        &     \\\hline
277.a.277.2        &     $277^{1}$        &$5^{1}$                         					&                        $2^{1}3^{3}31^{-2}$        &     0        &     	          6        &     \\
\hline
\end{tabular}
\label{gen2tab4}
\end{table}

\pagebreak

\begin{table}[h!]
\centering
\caption{$p=5,\,q=11,\,\psi(2)=\zeta_5,\,\p=(1-\zeta_5).$}
\scriptsize
\begin{tabular}{|c|c|c|c|c|c|c|} \hline

Curve Label                    &                   $N_A$                    &               $|A(\Q)_{\mathrm{tors}}|$                    &    $\cL(A,\triv)$                    &                   		$\cL(A,\psi)$                                     					&                                                             $\ord_\p$                              &$\Theta_C$ \\ \hline \hline
                                                                                                                                                                          
169.a.169.1                    &                   $13^{2}$                    &               $19^{1}$                                   			  &                                   $7^{1}11^{-2}19^{-1}$                    &                   $-(1+\zeta_{5}^2)$                                  				&                                                    1                    &                   \\\hline
196.a.21952.1                    &                   $2^{2}7^{2}$                    &             $2^{2}3^{2}$                    &                   $2^{2}3^{-1}11^{-2}$                    &                   $3\zeta_{5}^3$                                   				  &                                                 1                    &                  \\\hline
249.a.249.1                    &                   $3^{1}83^{1}$                    &            $2^{1}7^{1}$                                     	&                                          $2^{2}7^{-1}11^{-2}$                    &                    $-(1+\zeta_{5}+\zeta_{5}^2)$                    &                          1                    &                   \\
\hline

\end{tabular}
\label{gen2tab5}
\end{table}

\pagebreak

\begin{table}[h!]
\centering
\caption{$p=5,\,q=31,\,\psi(3)=\zeta_5,\,\p=(1-\zeta_5).$}
\scriptsize
\begin{tabular}{|c|c|c|c|c|c|c|} \hline
Curve Label        &       $N_A$        &    $|A(\Q)_{\mathrm{tors}}|$   &                         $\cL(A,\triv)$                         				&                               $\cL(A,\psi)$                      	&                            $\ord_\p$         &   $\Theta_C$ \\ \hline \hline
                                                                                                                   
169.a.169.1        &       $13^{2}$        & $19^{1}$        &    $2^{4}3^{1}19^{-1}31^{-2}$        &       $3(3-2\zeta_{5}+3\zeta_{5}^2)$        &                1        &              \\\hline
196.a.21952.1        &       $2^{2}7^{2}$   &      $2^{2}3^{2}$&   $2^{2}3^{1}31^{-2}$        &       $-3(8\zeta_{5}-5\zeta_{5}^2+8\zeta_{5}^3)$        &           1        &              \\\hline
249.a.249.1        &       $3^{1}83^{1}$           &     $2^{1}7^{1}$        &       $2^{5}7^{-1}31^{-2}$        &       $5+3\zeta_{5}^2+3\zeta_{5}^3$        &                    2        &              \\\hline
249.a.6723.1        &       $3^{1}83^{1}$          &     $2^{2}7^{1}$        &       $2^{5}7^{-1}31^{-2}$        &       $5+3\zeta_{5}^2+3\zeta_{5}^3$        &                    2        &               \\\hline
256.a.512.1        &       $2^{8}$           & $2^{2}{\color{blue}5^{1}}$  &            $2^{3}31^{-2}$        &       $1+\zeta_{5}+\zeta_{5}^2$                         				&                                             1        &               \\\hline
277.a.277.1        &       $277^{1}$               &        $3^{1}{\color{blue}5^{1}}$&   $2^{1}3^{1}31^{-2}$        &       $-(\zeta_{5}-3\zeta_{5}^2+\zeta_{5}^3)$  &                     1        &               \\\hline
277.a.277.2        &       $277^{1}$           & ${\color{blue}5^{1}}$ &     $2^{1}3^{3}31^{-2}$        &       $-9(\zeta_{5}-3\zeta_{5}^2+\zeta_{5}^3)$   &                   1        &               \\
\hline

\end{tabular}
\label{gen2tab6}
\end{table}

\pagebreak

\begin{table}[h!]
\centering
\caption{$p=7,\,q=29,\,\psi(2)=\zeta_7,\,\p=(1-\zeta_7).$}
\scriptsize
\begin{tabular}{|c|c|c|c|c|c|c|} \hline
Curve Label                                      &                   $N_A$                                      &            			 $|A(\Q)_{\mathrm{tors}}|$                                      &       $\cL(A,\triv)$                                      &                   $\cL(A,\psi)$                                                      								 &                                          	                   $\ord_\p$                                      &                   $\Theta_C$ \\ \hline \hline
                                                                                                                                                                                                               
169.a.169.1                                      &                   $13^{2}$                                      &                     $19^{1}$                                      &                            $2^{4}3^{1}19^{-1}29^{-2}$                                      &                   $-(1+\zeta_{7}+\zeta_{7}^4)$                                                       				&                                1                                      &                   \\\hline
196.a.21952.1                                      &                   $2^{2}7^{2}$                                      &                  $2^{2}3^{2}$                                      &                           $2^{2}3^{1}29^{-2}$                                      &                    $-12(1+\zeta_{7}+\zeta_{7}^2+\zeta_{7}^3+\zeta_{7}^4+\zeta_{7}^5)$                 &             1                                      &                   \\
\hline
\end{tabular}
\label{gen2tab7}
\end{table}

\pagebreak

\Addresses

\end{document}